\documentclass[12pt]{article}

\usepackage[T1]{fontenc}
\usepackage[T2A]{fontenc}
\usepackage[utf8]{inputenc}
\usepackage[russian,english]{babel}

\usepackage[mathscr]{eucal}
\usepackage{amsbsy}
\usepackage{amssymb}
\usepackage{amsmath}
\usepackage{amsthm}

\usepackage{latexsym}

\usepackage{stmaryrd}
\usepackage[dvips]{graphicx,epsfig}
\usepackage{color}
\newcommand{\HH}{\mathrm{H}}
\newcommand{\PP}{\mathbf{P}}

\newcommand{\ve}{\varepsilon}
\def\phi{\varphi}
\def\bbn{{\mathbb N}}
\def\bbr{{\mathbb R}}
\def\bbb{{\mathbb B}}

\def\bz{{\mathbb Z}} %Keistas Z apribrezimas

\def\Bl1{{\bf 1}}
\def\B2{{\bf 2}}
\def\B0{{\bf 0}}

\def\g{\gamma}
\def\l{\lambda}

\def\s{\sigma}

\def\=A8{\"o}

\def\co{\mathop{\hbox{\rm conv}}\nolimits}

\newcommand{\beq}{\begin{equation}}
\newcommand{\eeq}{\end{equation}}
\newcommand\beqn{\begin{displaymath}}  % no number
\newcommand\eeqn{\end{displaymath}}
 \newcommand{\sv}[1]{\left\lfloor #1 \right\rfloor}
\newcommand{\halmos}{\vspace{3mm} \hfill \mbox{$\Box$}\\[2mm]}
\theoremstyle{plain}
\newtheorem{teo}{Theorem}
\newtheorem{lem}[teo]{Lemma}

\newtheorem{prop}[teo]{Proposition}
\theoremstyle{definition}

\newtheorem{remark}[teo]{Remark}

\begin{document}

\title{More on the convergence of Gaussian convex hulls %\footnote{The first named author was supported by grant VIZIT-3-TYR-013 of   Lithuanian Research Council.}
\footnotemark[0]\footnotetext[0]{ \textit{Short title:} asymptotic
of convex hulls of Gaussian sequences }
\footnotemark[0]\footnotetext[0]{%
\textit{MSC 2000 subject classifications}. Primary 60G15, secondary
 60F15 .} \footnotemark[0] \footnotemark[0]\footnotetext[0]{ \textit{Corresponding
author:} Youri Davydov, Universit\'e de Lille, Laboratoire Paul
Painlev\'e,
 e-mail: youri.davydov@univ-lille.fr   }}

\author{ Youri Davydov$^{\text{\small 1}}$  and  Vygantas Paulauskas$^{\text{\small 2}}$ \\
{\small $^{\text{1}}$ Universit\'e de Lille, Laboratoire Paul
Painlev\'e}\\
{\small and Saint Petersburg state university}\\
 {\small $^{\text{2}}$ Vilnius University, Department
of Mathematics and Informatics}\\
 }

\maketitle

\begin{abstract}
A "law of large numbers" for consecutive convex hulls for weakly dependent  Gaussian sequences $\{X_n\}$, having the same marginal distribution, is  extended to the case when the sequence $\{X_n\}$ has a weak limit. Let  $\mathbb{B}$ be a separable Banach space with a conjugate space $\mathbb{B}^\ast$. Let $\{X_n\}$
be a centered  $\mathbb{B}$-valued Gaussian sequence satisfying two conditions: 1) $X_n \Rightarrow X\;\;$ and 2) For every $x^* \in \mathbb{B}^\ast$
$$
\lim_ {n,m, |n-m|\rightarrow \infty}E\langle X_n, x^*\rangle \langle X_m, x^*\rangle\;\; = \;\;0.
$$
Then with probability 1 the normalized convex hulls
$$
W_n = \frac{1}{(2\ln n)^{1/2}}\,\co \{\,X_1,\ldots,X_{n}\,\}
$$
converge in Hausdorff distance to the concentration ellipsoid of a limit Gaussian $\mathbb{B}$-valued random element $X.$
In addition,
some related questions are discussed.

\end{abstract}
%\vfill
%\eject
\section{Introduction and formulation of results}

 Let   $\mathbb{B}$ be a separable Banach space with a norm $||\cdot ||$ and let
$\mathbb{B}^\ast$ and  $\left\langle \cdot,\,\cdot\right\rangle$ denote its conjugate space and  the corresponding inner product, respectively.  For $A \in {\mathbb B}$ the notation $\co\{A\}$ is used for the closed convex hull of $A.$ If $X$ is a $\mathbb{B}$-valued centered Gaussian random element with a distribution ${\cal P}$ then  by $\HH$ we denote its reproducing kernel Hilbert space   and ${\cal E}$ will stand for the closed unit ball in $\HH$,  see, e.g., \cite{Ledoux}, p. 207. The set ${\cal E}$ is also called {\it concentration ellipsoid} of $X.$

Finally, we introduce the separable complete metric space
$\mathcal{K}_{\mathbb B}$ of all nonempty compact subsets of a Banach space $\mathbb B$ equipped with
the Hausdorff distance $d_{\mathbb B}$:
$$
d_{\mathbb B}(A,B) = \max\{\inf\{\,\epsilon \;|\;A \subset
B^\epsilon\},\;\; \inf\{\,\epsilon \;|\;B \;\subset
A^\epsilon\}\},
$$
where $A^\epsilon$ is the open $\epsilon$-neighborhood of $A$. Convergence
of compact sets in $\bbb$ always will be considered in this metric.

Investigation of the asymptotic behavior of convex hulls
$$
W_n = \co \{\,X_{1},\ldots,X_n\}
$$
of multivariate Gaussian random variables is an important part of Extreme Value Theory and has various applications, see for example, \cite{Randon} and reference list, containing 160 items, in it.
In 1988 Goodman \cite{Goodman} proved a fundamental result that the normalized set  $\{\,X_{1},\ldots,X_n\}$ of independent and identically distributed $B$-valued centered Gaussian random elements with a distribution ${\cal P}$ is approaching the concentration ellipsoid of ${\cal P}$  as $n$ grows to  infinity. From this result one can immediately derive that a.s.
\begin{equation}\label{conv}
\frac{1}{b(n)}W_n \stackrel{\mathcal{K}_{\mathbb B}}{\rightarrow} {\cal E}, \quad {\rm as} \ n\to \infty,
\end{equation}
where $b(t) = \sqrt{2\ln (t)}, \ t>e$. Moreover, the rate of convergence in this relation is of the order $o(b(n)^{-1}).$

Later the convergence of the type  (\ref{conv}) was proved  first for stationary $d$-dimensional weakly dependent Gaussian sequences in \cite{DavDom}, and then  the similar result was proved for $\mathbb{B}$-valued  Gaussian random fields on $\bbr^m$ or $\bz^m$ in \cite{PaulDav1}. Although at the introduction of the paper \cite{PaulDav1} it was said that only the case $m>1$ is considered, inspection of the proof of the main result -Theorem 1.1 in \cite{PaulDav1} - shows that the result holds for $m=1$, too. In particular, this result states that if a $\mathbb{B}$-valued centered Gaussian  sequence  $\{X_k, k\in \bbn\}$ has the same marginal distribution ${\cal P}$ and satisfies  the following condition
\begin{equation}\label{cond2}
E\left\langle X_n,\,x^\ast\right\rangle\left\langle
X_m,\,x^\ast\right\rangle \rightarrow 0,  \quad {\rm  for \ all}\ x^\ast \in \mathbb{B}^\ast \ \ {\rm as} \ n, m, |n-m|\to \infty.
\end{equation}
then  (\ref{conv}) holds.

In the paper we show that the condition of equality of marginal distributions  can be essentially relaxed substituting it by the weak convergence of the sequence $\{X_n\},$ and for weak convergence we use the sign $\;\Rightarrow.$

\begin{teo}\label{thm1} Suppose that a centered Gaussian sequence of $\mathbb{B}$-valued random elements $\{X_k, k\in \bbn\}$ satisfies $(\ref{cond2})$ and the following condition:
\begin{equation}\label{cond1}
 X_n \Rightarrow X.
\end{equation}
Then a.s.
\begin{equation}\label{convWn}
\frac{1}{b(n)}W_n \stackrel{\mathcal{K}_{\mathbb B}}{\rightarrow} {\cal E}, \quad {\rm as} \ n\to \infty,
\end{equation}
where ${\cal E}$ is concentration ellipsoid of $X.$
\end{teo}

Since the proof of this theorem will be carried in two steps, and in the first step we consider the case $\mathbb B=\bbr$, we look more closely what is the meaning of the result in  this particular case. Let $N(0, \s^2)$ stand for a Gaussian random variable with mean zero and variance $\s^2$, and $\{X_k\},\ k\in \bz_+,$ is a sequence of $N(0, \s_k^2)$ random variables. Without loss of generality we can assume that $X$ in (\ref{cond1}) is $N(0,1)$. We have
$$
W_n=[\min\{X_1, X_2, \dots , X_n\}, \max\{X_1, X_2, \dots , X_n\}]
$$
 and ${\cal E}=[-1, 1]$. If the covariance function $\rho (m, n):=EX_mX_n \to 0$ as
$ n, m, |n-m|\to \infty$, then we have the relation (\ref{convWn}). It is clear, that if the dependence between elements $X_k$ of the sequence is stronger, the sequence of their convex hulls  is more concentrated. One can consider the extreme case, when $X_k\equiv X_0$ for all $k\ge 1$, then $W_n=\{X_0\}$ is one point and $\lim_{n\to \infty}(g(n))^{-1} W_n =\{0\}$ for any sequence $g(n)\to \infty$. The following example gives us additional information in this question.
\begin{remark}\label{examp}Let us consider the sequence of i.i.d. $N(0, 1)$ random variables $\{\xi_j\}, \ j\ge 1,$ and let $S_k=k^{-1/2}\sum_{j=1}^k \xi_j$. Taking $X_k=S_k,$ we are in the setting of Theorem \ref{thm1}, but the condition (\ref{cond2}) is not satisfied, since, if $n=m+k, k>0$, then
\begin{equation}\label{rhomn}
\rho (m, n)=\frac{ES_mS_{m+k}}{(m(m+k))^{1/2}}=\frac{m}{(m(m+k))^{1/2}}=\left (1+\frac{k}{m} \right )^{-1/2}.
\end{equation}
Thus, in order to get $\rho(m, n)\to 0$, it is not sufficient to require $m\to \infty, \ k\to \infty$, but stronger condition is required $m\to \infty, \ k/m \to \infty.$ On the other hand, $S_k$ is a sum of i.i.d.random variables, therefore, denoting $c(t)=(2\ln\ln t)^{1/2}, \ t> e,$ the classical LIL gives us that the cluster set  for the sequence $\{X_n/c(n)\}$ is $[-1, 1]$, while for the sequence $\{X_n/b(n)\}$ with probability one limit is zero. We shall prove that for this example   we have the following result:
\begin{prop}
 With probability one
\begin{equation}\label{convWnR1}
\frac{1}{c(n)}W_n \stackrel{\mathcal{K}_{\mathbb R}}{\rightarrow} [-1, 1], \quad {\rm as} \ n\to \infty,
\end{equation}
where $W_n=[-V(n), V(n)]$, and $V(n):=\max \{X_1, X_2, \dots , X_n\}$
\end{prop}
\end{remark}

In connection with this example it is possible to formulate the following problem.

 Suppose that a sequence $\{X_n\}$ has standard normal marginal distributions and covariance function $\rho (m, n )$. For which functions $g(n)$ and under what conditions for covariance function $\rho$ we can get  the relation (\ref{convWnR1}) with function $g$ instead of $c$?

  This Proposition and Theorem \ref{thm1} give us two examples of such functions $g$. What other normalizing functions are possible in relation (\ref{convWnR1})?

Let us make three final remarks.
\begin{remark} As in \cite{PaulDav1}, having (\ref{convWn}), we can get some information on asymptotic behavior of $Ef(W_n)$ for some functions defined on $\mathcal{K}_{\mathbb B}$. In case ${\mathbb B}=\bbr^m$ typical examples of such functions are diameter, volume or surface measure.
\end{remark}

\begin{remark} From the proof of Theorem \ref{thm1} we can extract some information about the rate of convergence in (\ref{convWn}). Namely, in the proof we have the equality
$$
\mathbb{P}\left \{d_{\mathbb B}(X, b(n){\cal E}) >\varepsilon\right \}=\mathbb{P}\{\exp{\frac{1}{2}\psi_\ve}\ge n\},
$$
and since these probabilities are monotonically non-increasing and $\sum_n\mathbb{P}\{\exp{\psi_\ve/2}\ge n\}=E\{\exp{\psi_\ve/2}\}<\infty $ we get
$$
 \mathbb{P}\left \{d\left (\frac{X_n}{b_n},{\cal E}\right ) \,> \,\varepsilon \right \} = o(n^{-1}).
$$
Let us note that this result cannot be compared with the result from \cite{Goodman}, where it is proved that with probability $1$
$$
d_{\mathbb B}\left(\frac{X_n}{b_n},{\cal E}\right) = o(b_n^{-1}).
$$
\end{remark}

\begin{remark} In Theorem \ref{thm1} and in previous results for Gaussian sequences limit set of convex hulls was ellipsoid of some Gaussian measure. If we dismiss the condition of weak convergence of Gaussian sequence  $\{X_k\},$ the limit set may exist, but not necessarily will be an ellipsoid. For example, the following statement holds.

\begin{prop}\label{prop5}
Let $\mathbb{B}$ be a  separable Banach space. Let $V\subset \mathbb{B}$ be  a central symmetric polytope,
$V = \co\{a_k, -a_k,\, k = 1,\ldots,m\}.$
%convex hull of the union of  $m$  non-collinear segments $[-a_k,a_k].$
Then there exists a sequence of independent Gaussian vectors $\{X_k\}$ such that a.s.
$$
\frac{1}{b(n)}W_n\;\rightarrow \;V.
$$
\end{prop}

%\vspace{5pt}

\end{remark}

\section{Proofs }

{\it Proof of Theorem $\ref{thm1}$}.
 As it was mentioned above, the proof will be carried in two steps, and in the first step we consider the case $\mathbb B=\bbr$.
\vspace{5pt}

I. Without loss of generality we can suppose that $X$ has a standard Gaussian distribution; then ${\cal E}=[-1, 1]$.

The condition (\ref{cond2}) transforms now in
$$
EX_nX_m \rightarrow 0\;\;\; \mathrm{as}\;\;\;n,m,|n-m| \rightarrow \infty.
$$
From weak convergence of $X_n$ to $X$ follows that $\sigma_n^2 :=EX_n^2 \rightarrow 1$ as $n\rightarrow \infty.$

For r.v. $Y_n = X_n/\sigma_n$ the conditions (\ref{cond1}) and (\ref{cond2}) are fulfilled. Since $Y_n$ are identically distributed,
setting $U_n = \co\{Y_1,\ldots,Y_n\},$
by Theorem 1.1. from \cite{PaulDav1}
\begin{equation}\label{convR1}
\frac{1}{b(n)}U_n \;\;\stackrel{\mathcal{K}_{1}}{\longrightarrow}\;\;\; [-1,1],\;\; a.s. \quad {\mathrm as}\, \ n\to \infty.
\end{equation}
Let us show that a.s.
\begin{equation}\label{convR2}
\Delta_n:= d_{{\bbr}^1}\left(\frac{1}{b(n)}W_n,\;\frac{1}{b(n)}U_n\right)
\rightarrow 0.
\end{equation}
We have
\begin{equation}\label{convR3}
\Delta_n \leq \frac{1}{b(n)}\max_{k\leq n}\{|X_k-Y_k|\}=
\frac{1}{b(n)}\max_{k\leq n}\{|Y_k||1-\sigma_k|\}.
\end{equation}
Let $Z_n=\max_{k\leq n}\{|Y_k|\};\;\;\;M_n=\max_{k\leq n}\{|1-\sigma_k|\}.$ For $\varepsilon>0$ find $n_0$ such that $\sup_{k>n_0}|1-\sigma_k|< \varepsilon.$
Then for $n\geq n_0$ by (\ref{convR3})
$$
\Delta_n \leq \frac{1}{b(n)}\max_{n_0\leq k\leq n}\{|Y_k|\}\varepsilon
+ \frac{M_{n_0}Z_{n_0}}{b(n)}.
$$
It follows from (\ref{convR1}) that $\limsup_n \Delta_n \leq \varepsilon.$ Hence we get (\ref{convR2}) and (\ref{convWn}) is proved for $\mathbb{B}=\mathbb{R}.$
\vspace{5pt}

II. General case. We shall show that with the probability 1 the sequence
$\{b(n)^{-1}W_n\}$ is relatively compact in ${\mathcal{K}_{\mathbb B}}.$ Due to Lemmas 2.2 and 2.3 from \cite{PaulDav1} it is sufficient to prove that there exists a compact set $K$ such that for
every $\varepsilon>0$ with probability 1 for all sufficiently large $n$ we have the following inclusion
$\{b(n)^{-1}W_n\}\subset K^\varepsilon.$
We take $K = {\cal E}.$ Since the space ${\mathbb B}$ is fixed, instead of $d_{\mathbb B}$ we shall write simply $d.$ It is clear that this inclusion will  follow from the following relation:  for every $\varepsilon >0$ a.s.
\begin{equation}\label{comp1}
\limsup_n d(X_n, (1+\varepsilon)b(n){\cal E}) = 0.
\end{equation}

By Skorokhod representation theorem we can suppose that
$X_n \rightarrow  X$ a.s.

Let $\sigma^2 = E\|X\|^2,\;\sigma_{max}^2 = \sup_n E\|X_n\|^2.$ It follows from Fernique's theo\-rem about
integrability of exponential moments (see \cite{Fernique}) that for every
$\gamma ,\;0<\gamma<(2\sigma^2)^{-1},$
$$%\begin{equation}\label{comp2}
E\exp\{\gamma\|X\|^2\} \;<\;\infty.
$$%\end{equation}
Moreover, from the proof of Fernique's theorem one can deduce that
for every $\gamma ,\;0<\gamma<(2\sigma_{max}^2)^{-1},$
$$%\begin{equation}\label{comp3}
\limsup_n E\exp\{\gamma\|X_n\|^2\} \;<\;\infty.
$$%\end{equation}
It means, in particular, that for every $\gamma ,\;0<\gamma<(2\sigma_{max}^2)^{-1},$ the family $\exp\{\gamma\|X_n\|^2\}$ is uniformly
integrable. Therefore
$$
\delta_n^2:= E\|X_n -X\|^2 \;\rightarrow \;0\;\;{\mathrm {and}}\;\;
E\exp\{\gamma\|X_n-X\|^2\}\;\rightarrow \;\;1.
$$
For $\varepsilon >0$ let $n_\varepsilon$ be such that
$\sigma_\varepsilon^2:= \sup_{n>n_\varepsilon}\delta_n^2\,<\,\varepsilon^2.$

 We have
 $$
 \mathbb{P}\{d(X_n, (1+\varepsilon)b(n){\cal E}) \,> \,\varepsilon\}\,=\,A_n+B_n,
 $$
where
$$
A_n = \mathbb{P}\{d(X_n, (1+\varepsilon)b(n){\cal E}) \,> \,\varepsilon,\;\|X_n-X\|<\varepsilon b(n)\},
$$
$$
B_n =\mathbb{P}\{d(X_n, (1+\varepsilon)b(n){\cal E}) \,> \,\varepsilon,\;\|X_n-X\|\geq\varepsilon b(n)\}.
$$
Evidently
$$
A_n \leq \mathbb{P}\{d(X, b(n){\cal E}) >\varepsilon\}.
$$
Now we formulate Talagrand's lemma \cite{Talagrand} as it is formulated   in \cite{Goodman}, see Lemma 3.1 therein.

\begin{lem}\label{lem}
Let $X$ be a $\mathbb{B}$-valued centered Gaussian random element with a  concentration ellipsoid ${\cal E}$. Then for any $\varepsilon>0$ there is a random variable
$\psi_\ve$ such that
$$
E\{\exp{\{\frac{1}{2}\psi_\ve}\}\}<\infty
$$
 and for all $\l>0$
$$
\mathbb{P}\{d(X, \l{\cal E}) \le\varepsilon\}=\mathbb{P}\{\psi_\ve<\l^2 \}.
$$
\end{lem}
We apply this Lemma taking $\l=b(n)$ and obtain
$$
A_n\leq \mathbb{P}\{d(X, b(n){\cal E}) >\varepsilon\}=\mathbb{P}\{\psi_\ve\ge 2\ln n \}=\mathbb{P}\{\exp{\frac{1}{2}\psi_\ve}\ge n\}
$$
Hence $\sum_nA_n < \infty.$

For $0<\g< (2\sigma_\varepsilon^2)^{-1}$ we have $E\exp\{\g\|X_n-X\|^2\}<\infty$ for each $n>n_\ve$, therefore, denoting
$$
L_\ve(a) = \sup_{n>n_\ve}\{E\exp\{a\|X_n-X\|^2\}\},
$$
 for
$a\in \left(\frac{1}{2\ve^2},\frac{1}{2\sigma_\ve^2} \right)$ and $n \ge n_\ve$
we  apply once more Fernique's theorem and get
$$
B_n \leq \mathbb{P}\{\|X_n-X\|>\ve b(n)\} \leq
\frac{E\exp\{a\|X_n-X\|^2\}}{n^{2a\ve^2}} \leq \frac{L_\ve(a)}{n^{2a\ve^2}}.
$$

Since $2a\ve^2>1,$ this estimate gives the convergence of the series
$\sum_nB_n.$

 Therefore we see that for every $\ve>0$
 $$
 \sum_n \mathbb{P}\{d(X_n, (1+\varepsilon)b(n){\cal E}) \,> \,\varepsilon\} <\infty.
 $$

Then the Borel-Cantelli lemma  gives us (\ref{comp1}),
which shows that for every $\delta >0$ with probability 1 for all sufficiently large $n$
$$
\frac{1}{b(n)}W_n \subset {\cal E}^\delta.
$$
This proves the relative compactness of \{$b(n)^{-1}W_n\}.$

It follows from Lemma 2.7 \cite{PaulDav1} that now it is sufficient to prove the convergence for every $\theta \in S_1^*(0):=\{x^* \in {\mathbb B}^*: ||x^*||=1 \}$
\begin{equation}\label{convsupp}
M_n(\theta) \stackrel{a.s.}{\rightarrow} M_{\cal E}(\theta),\;\;\;n\rightarrow
\infty,
\end{equation}
where $M_n, M_{\cal E}$ are support functions for $b(n)^{-1}W_n$ and $  {\cal E}$, respectively. We recall that a function $M_A(\theta)$, defined by the relation
$$
M_A(\theta):=\sup_{x\in A}\langle x,\theta\rangle, \quad A\in \mathcal{K}_{\mathbb B}, \ \theta \in S_1^*(0),
$$
is called the support function of a set $A\in \mathcal{K}_{\mathbb B}.$
Since
$$
M_n(\theta) = \frac{1}{b(n)}\max_{k\leq n}{\{\langle X_k,\theta\rangle\}},
$$
the convergence (\ref{convsupp}) follows from the first part of the proof. \halmos

{\it Proof of Proposition $\ref{examp}$}. Let us denote
$$
c_1:=\liminf_n \frac{V(n)}{c(n)}.
$$
Let $l\ge 2$ be a fixed integer, then we have $c(l^m)\sim b(m)$, as $m\to \infty.$  $V(n)$ is non-decreasing, therefore, for $n\in [l^k, l^{k+1}]$, we have
$$
\frac{V(n)}{c(n)}\ge \frac{c(l^k)}{c(n)}\frac{V(l^k)}{c(l^k)}.
$$
Since
$$
\frac{c(l^k)}{c(n)}\ge \frac{c(l^k)}{c(l^{k+1})}\to 1, \ \ {\rm as} \  n\to \infty,
$$
we get
$$
c_1\ge \liminf_m \left\{ \frac{V(l^m)}{c(l^m)}\right\}\ge
\liminf_m \left\{\frac{\max\{X_{l}, X_{l^2}, \dots ,  X_{l^m}\}}{c(l^m)}\right\}.
$$
In our example we have (\ref{rhomn}), therefore $r:=\sup_{i\ne j}EX_{l^i}X_{l^j}=l^{-1/2}$. Due to Lemma 2.5 from \cite{PaulDav1} we get
$$
\liminf_m \left\{\frac{X_{l}, X_{l^2}, \dots ,
 X_{l^m}}{c(l^m)}\right\}\ge \sqrt {1-r}=\sqrt {1-l^{-1/2}}.
$$
This quantity can be made close to $1$ if we choose $l$ sufficiently large, therefore with probability one we have
\begin{equation}\label{estbel}
c_1 \ge 1.
\end{equation}
In order to get the estimate from above for
$$
c_2 :=\limsup_n \left\{\frac{V(n)}{c(n)}\right\}
$$
we shall need the following lemmas.
\begin{lem}\label{lem1} Suppose that a sequence of random variables $\{Y_k\}$ satisfies the following condition: for all $\g<(2\s^2)^{-1},\;
\s >0, $
\begin{equation}\label{EgY}
\sup_n E\exp \{\g Y_n^2\}<\infty.
\end{equation}
Then
$$
\limsup_n \left\{\frac{\max_{k\le n}\{Y_k\}}{b(n)}\right\}\le \s.
$$
\end{lem}
The proof of this lemma coincides with the proof of Lemma 1 in \cite{Davydov}, despite of the fact that in this paper the variables $\{Y_k\}$ was i.i.d. It turns out that independence is not used at all and condition of identical distributions of $Y_k$ can be replaced by condition (\ref{EgY}).

\begin{lem}\label{lem2} $(\cite{Petrov},$ Thm. $2.2)$ Let $\{\xi_k\}, \ k\ge 1,$ be independent symmetric random variables, $S_n=\sum_{k=1}^n \xi_k$. Then for every $x\ge 0$
$$
\PP\left (\max_{1\le k\le n}|S_k|\ge x\right )\le 2\PP(|S_n|\ge x).
$$
\end{lem}
\begin{lem}\label{lem3} Let $X$ and $Y$ be two non-negative random variables with distribution functions $F$ and $G$. If for some $a\ge 1$ and for all $x\ge 0$ we have $1-F(x)\le a (1-G(x)).$
Then for every $\g>0$
$$
E\exp \{\g X^2\}\le aE\exp \{\g Y^2\}.
$$
\end{lem}
Elementary proof of this statement follows from equalities
 $$
 E\exp \{\g X^2\}=\sum_{k=0}^\infty\g^k(k!)^{-1}E(X^{2k}), \; E(X^{2k})=\int_0^\infty (1-F(x))d(x^{2k}).
 $$

Let us  fix a non-integer (we have in mind  that we shall choose $a$ close to $1$) number $a>1$ and let us denote
$$\Delta_j=\{i\in N: a^j\le i\le a^{j+1}\},\;\ Y_j=\max_{i\in \Delta_j}X_i.
$$
 As in the case of lower bound  we can prove that
$$
c_2 \le \limsup_m \left\{\frac{V(\sv{a^m})}{c(\sv{a^m})}\right\}.
$$
Note that
$$
\Delta_j=\{i\in N: \sv{a^j}+1\le i\le \sv{a^{j+1}}\}
$$
and $V(\sv{a^m})=\max_{0\le j\le m-1}Y_j$,  therefore
$$
|Y_j|\le \max_{i\in \Delta_j}\{|X_i|\}=\max_{i\in \Delta_j}\left\{|\frac{S_i}{\sqrt i}|\right\}\le \frac{1}{\sqrt {\sv{a^j}+1}}\max_{i\in \Delta_j}\{|S_i|\}\le \frac{1}{\sqrt {a^j}}\max_{i\in [\sv{a^j}+1, \sv{a^{j+1}}]}\{|S_i|\}.
$$
Applying Lemma \ref{lem2}, we have
$$
\PP\{\max_{i\in \Delta_j}|S_i|\ge x\}\le 2\PP\{|S_{\sv{a^{j+1}}}|\ge x\},
$$
whence
$$
\PP\{|Y_j|\ge x\}\le 2\PP\left \{\Big |\frac{S_{\sv{a^{j+1}}}}{\sqrt {a^j}}\Big |\ge x \right \}.
$$
Let  $\xi(j, a)=(a^j)^{-1}S_{\sv{a^{j+1}}}$. It is easy to see that $\xi(j, a)$ has distribution $N(0, \s^2(j,a))$ with
\begin{equation}\label{estsigma}
\s^2(j,a)=\frac{\sv{a^{j+1}}}{a^j}\to a, \ {\rm as} \ j\to \infty, \quad {\rm and} \ \ \s^2(j,a)\le a.
\end{equation}
Applying Lemma \ref{lem3} with $\g<1/2a,$ we get
$$
E\exp \{\g Y_j^2\}\le 2E\exp \{\g \xi(j, a)^2\}.
$$
Due to (\ref{estsigma}) we have $\sup_j E\exp \{\g \xi(j, a)^2\}:=C(a)<\infty$, therefore, using Lemma \ref{lem1} with $\s^2=a$ and recalling that  $c(\sv{a^m})\sim b(m)$, as $m\to \infty$, we get that with probability 1
$$
c_2\le \limsup_m \frac{1}{b(m)}V(\sv{a^m})\le {\sqrt a}.
$$
Since the last estimate holds  for any $a>1$, we get that with probability 1
\begin{equation}\label{estabove}
c_2 \le 1.
\end{equation}
Estimates (\ref{estbel}) and (\ref{estabove}) prove (\ref{convWnR1}).
\halmos
\vspace{5pt}

{\it Proof of Proposition $\ref{prop5}$.}
Let $\bbn=\cup_{k=1}^m T_k$, where the sets $T_k, \ k=1, \dots , m,$ are disjoint and have positive densities $p_k$. Let  $\{X_n\}$ be a sequence of independent random vectors such that for each $k$ and
$ j\in T_k,\;\;X_j$ has  Gaussian distribution concentrated on the line  $\{ta_k, t\in R^1\}$ with zero mean and variance $\sigma_k^2 =
\parallel a_k\parallel.$
We denote  $W_n^{(k)} = \co\{X_j\;j\leq n,\; j\in T_k\}.$
Since for any $p>0,$
$$
\lim_n \frac{b(np)}{b(n)} = 1,
$$
 Theorem \ref{thm1} implies that a.s. for any $k=1,\ldots,m,$
$$
\frac{1}{b(n)}W_n^{(k)}\;\rightarrow \;\co\{-a_k,a_k\}.
$$
Clearly, we have  $W_n = \co\{W_n^{(1)},\ldots,W_n^{(m)}\},$
therefore a.s.
$$
\frac{1}{b(n)}W_n\;\rightarrow \;\co\{\cup_{k=1}^m \co\{-a_k,a_k\}\}\,=\,V.
$$
\halmos

\section*{References}

\bibliographystyle{plain}

\end{document}